\documentclass{amsart}
\usepackage{amscd,amssymb,hyperref}
\theoremstyle{plain}
\newtheorem{prop}[subsection]{Proposition}
\newtheorem{thm}[subsection]{Theorem}
\newtheorem{lem}[subsection]{Lemma}
\newtheorem{cor}[subsection]{Corollary}

\theoremstyle{remark}
\newtheorem{rem}[subsection]{Remark}

\theoremstyle{definition}
\newtheorem{defn}[subsection]{Definition}
\newtheorem{exm}[subsection]{Example}

\numberwithin{equation}{section}

\renewcommand{\b}[1]{\mathbf{#1}}

\newcommand{\A}{{\mathcal A}}
\newcommand{\B}{{\mathcal B}}
\newcommand{\cQ}{{\mathcal Q}}
\newcommand{\E}{{\mathsf E}}
\newcommand{\I}{{\mathsf I}}
\newcommand{\sA}{{\mathsf A}}
\newcommand{\cL}{{\mathcal L}}
\newcommand{\sL}{{\mathsf L}}
\newcommand{\bE}{{\mathbb E}}
\newcommand{\bF}{{\mathbb F}}
\newcommand{\bH}{{\mathbb H}}
\newcommand{\Z}{{\mathbb Z}}
\newcommand{\C}{{\mathbb C}}
\newcommand{\F}{{\mathbb F}}

\renewcommand{\a}{{\alpha }}
\renewcommand{\c}{{\gamma }}

\DeclareMathOperator{\rank}{rank}
\DeclareMathOperator{\id}{id}
\DeclareMathOperator{\ad}{ad}
\DeclareMathOperator{\Der}{Der}
\DeclareMathOperator{\codim}{codim}
\DeclareMathOperator{\Tors}{Torsion}
\DeclareMathOperator{\Prim}{Prim}
\DeclareMathOperator{\Symp}{Sp}

\begin{document}

\title[fiber-type arrangements]
{Lie Algebras Associated to Fiber-Type Arrangements}

\author[D.~Cohen]{Daniel C.~Cohen$^\flat$}
\address{Department of Mathematics, Louisiana State University, Baton
Rouge, LA 70803}
\email{\href{mailto:cohen@math.lsu.edu}{cohen@math.lsu.edu}}
\urladdr{\href{http://www.math.lsu.edu/~cohen/}
{http://www.math.lsu.edu/\~{}cohen}}
\thanks{{$^\flat$}Partially supported by grants LEQSF(1996-99)-RD-A-04
and LEQSF(1999-2002)-RD-A-01 from the Louisiana Board of Regents, and
by grant MDA904-00-1-0038 from the National Security Agency}

\author[F.~Cohen]{Frederick R.~Cohen$^\natural$}
\address{Department of Mathematics, University of Rochester,
Rochester, NY 14627}
\email{\href{mailto:cohf@math.rochester.edu}{cohf@math.rochester.edu}}
\urladdr{\href{http://www.math.rochester.edu/u/cohf/}
{http://www.math.rochester.edu/u/cohf}}
\thanks{{$^\natural$}Partially supported by the National Science
Foundation}

\author[M.~Xicot\'encatl]{Miguel Xicot\'encatl$^\sharp$}
\address{Depto. de Mathem\'aticas, Cinvestav del IPN,
Mexico City \phantom{123456789012345678901234567}
\phantom{23}
Max-Plank-Institut f\"ur Mathematik,
P.O.~Box 7280, D-53072 Bonn, Germany}
\email{\href{mailto:xico@mpim-bonn.mpg.de}{xico@mpim-bonn.mpg.de}}
\thanks{{$^\sharp$}Partially supported by the CINVESTAV of the IPN
and by the MPI f\"ur Mathematik}



\begin{abstract}
    Given a hyperplane arrangement in a complex vector space of
    dimension $\ell$, there is a natural associated arrangement of
    codimension $k$ subspaces in a complex vector space of dimension
    $k \ell$.  Topological invariants of the complement
    of this subspace arrangement are related to those of the
    complement of the original hyperplane arrangement.  In particular,
    if the hyperplane arrangement is fiber-type, then, apart from
    grading, the Lie algebra obtained from the descending central
    series for the fundamental group of the complement of the
    hyperplane arrangement is isomorphic to the Lie algebra of
    primitive elements in the homology of the loop space for the
    complement of the associated subspace arrangement.  Furthermore,
    this last Lie algebra is given by the homotopy groups modulo
    torsion of the loop space of the complement of the subspace
    arrangement.  Looping further yields an associated Poisson
    algebra, and generalizations of the ``universal infinitesimal
    Poisson braid relations.''
\end{abstract}


\maketitle

\section{Introduction} \label{sec:intro}

Two classical constructions of interest in group theory and topology
are:
\begin{enumerate}
\item[(i)] The Lie algebra arising from the filtration quotients
associated to the descending central series of a discrete group $G$;
and

\smallskip

\item[(ii)] The Lie algebra of primitive elements in the singular
homology of the loop space of a space $X$, for certain topological
spaces $X$.
\end{enumerate}
The purpose of this article is to illustrate that these two a priori
unrelated Lie algebras are, in fact, isomorphic in certain natural
cases.  This work is motivated by recent results relating the Lie
algebras of (i) and (ii) arising in the context of classical
configuration spaces, and resolves a conjecture of the second two
authors concerning the generalization of these results to spaces
arising from certain hyperplane arrangements.

The main result here asserts that the Lie algebra associated to the
fundamental group $G$ of the complement of a fiber-type hyperplane
arrangement is, up to regrading, isomorphic to the Lie algebra of
primitive elements in the homology of the loop space of the complement
of a higher dimensional analogue of the arrangement.  The main theorem
is, in fact, stronger.  The Samelson product for the loop space gives
rise to a graded Lie algebra given by the homotopy groups modulo
torsion.  This Lie algebra is, again up to regrading, also isomorphic
to the Lie algebra associated to the descending central series
quotients of~$G$.  In addition, after looping further, there are
natural related Poisson algebras arising from the homology of
associated iterated loop spaces.

Given a discrete group $G$, let $G_n$ be the $n$-th stage of the
descending central series, defined inductively by $G_1=G$, and
$G_{n+1}=[G_n,G]$ for $n\ge 1$, and let $E^n_0(G)=G_n/G_{n+1}$ be the
$n$-th associated quotient.  Let $E^*_0(G)=\bigoplus_{n\ge 1}
E^n_0(G)$ be the Lie algebra obtained from the descending central
series of $G$, with Lie algebra structure induced by the commutator
map $G\times G \to G$, $(x,y) \mapsto xyx^{-1}y^{-1}$.  For each
positive integer $k$, use the ungraded Lie algebra $E^*_0(G)$ to
define a related graded Lie algebra as follows.

\begin{defn} \label{defn:EGk}
For a group $G$, let $E^*_0(G)_k$ be the graded Lie algebra given by
\[
E^q_{0}(G)_k =
\begin{cases}
E^n_0(G) & \text{if $q=2nk$,}\\
0           & \text{otherwise,}
\end{cases}
\]
with Lie bracket structure induced by that of the Lie algebra
$E^*_0(G)$ obtained from the descending central series of $G$ in the
obvious manner.
\end{defn}

A theorem relating the Lie algebras of (i) and (ii) above is described
next.  Let $P_n$ be the Artin pure braid group, the fundamental group
of the configuration space $F(\C,n)$.  The results on configuration
spaces alluded to above, due to Fadell and Husseini \cite{FH} and
Cohen and Gitler \cite{CG}, may be summarized as follows.

\begin{thm} \label{thm:config}
For $k \ge 1$, the homology of the loop space of the configuration
space $F(\C^{k+1},n)$ is isomorphic to the universal enveloping
algebra of the graded Lie algebra $E^*_0(P_n)_k$.  Moreover,
\begin{enumerate}
\item[(a)] The image of the Hurewicz homomorphism
\[
\pi_*(\Omega F(\C^{k+1},n)) \to H_*(\Omega F(\C^{k+1},n);\Z)
\]
is isomorphic to $E^*_0(P_n)_k$; and

\smallskip

\item[(b)] The Hurewicz homomorphism induces isomorphisms of graded
Lie algebras
\[
\pi_*(\Omega F(\C^{k+1},n))/\Tors \to \Prim H_*(\Omega
F(\C^{k+1},n);\Z) \cong E^*_0(P_n)_k,
\]
where $\Prim\bullet$ denotes the module of primitive elements, and the
Lie algebra structure of the source is induced by the classical
Samelson product.

\end{enumerate}
\end{thm}

The Lie algebra arising in the above theorem is the ``universal
Yang-Baxter Lie algebra'' $\cL(n)$, the quotient of the free Lie
algebra on a free abelian group of rank $\binom{n}{2}$ by relations
recorded in \eqref{eq:braidrels} below, and known variously as the
``infinitesimal pure braid relations'' or the ``horizontal four-term
relations and framing independence.''  Furthermore, the homology of an
iterated loop space of configuration space, $\Omega^{q}F(\C^{k+1},n)$
for $q>1$, admits the structure of a graded Poisson algebra, see
\cite{CLM}.  The associated relations are called the ``universal
infinitesimal Poisson braid relations.''  This Poisson algebra
structure on $H_{*}(\Omega^{q}F(\C^{k+1},n))$ has recently been used
in the context of algebraic groups by Lehrer and Segal \cite{LS}.

The universal Yang-Baxter Lie algebra, and infinitesimal pure braid
relations, arise in a number of contexts.  These include the
classification of pure braids by Vassiliev invariants, see Kohno
\cite{K3}, and the Knizhnik-Zamolodchikov differential equations from
conformal field theory, where the relations appear as integrability
conditions on the associated Gauss-Manin connection, see Varchenko
\cite{Va}.  Moreover, any finite dimensional representation of the Lie
algebra $\cL(n)$ induces a representation of the pure braid group
$P_n$ on the same vector space, see Kapovich and Millson~\cite{KM}.

An important ingredient in the proof of Theorem~\ref{thm:config} is a
classical result of Fadell and Neuwirth \cite{FN} which shows that
configuration spaces admit iterated bundle structure.  Similar results
are known to hold for certain orbit configuration spaces
\cite{Xi,DCmono,CX}, which admit analogous bundle structure, and are
described in more detail below.  All of these spaces fit in the
following general framework.

For each natural number $\ell$, let $X_\ell$ be a functor from
Euclidean spaces, with morphisms restricted to endomorphisms, to
topological spaces.  For a Euclidean space $\bE$, let $\cQ_\ell(\bE)$
be a discrete subset of $\bE$ of fixed (possible infinite) cardinality
depending on $\ell$.  Assume that there are natural transformations
$X_{\ell}(\bE) \to X_{\ell-1}(\bE)$ which satisfy the following conditions.

\smallskip

\begin{enumerate}

\item The space $X_1(\bE) = \bE\setminus\cQ_1(\bE)$ is the complement
of a discrete subset of $\bE$.

\smallskip

\item The map $X_\ell(\bE) \to X_{\ell-1}(\bE)$ is a fiber bundle
projection, with fiber $\bE \setminus \cQ_\ell(\bE)$.

\smallskip

\item Each bundle $X_\ell(\bE) \to X_{\ell-1}(\bE)$ admits a
cross-section.

\smallskip

\item If $\bE\cong\C$, the fundamental group of $X_{\ell-1}(\bE)$ acts
trivially on the homology of the fiber $\bE\setminus\cQ_\ell(\bE)$.

\end{enumerate}
The prototypical examples are given by the configuration spaces
$X_{\ell}(\bE)=F(\bE,\ell)$, where $\bE=\C^{k}$.  Further examples are
given below.

It seems likely that for many choices of $X_\ell$, the Lie algebras
associated to $X_{\ell}(\bE)$ as $\bE$ varies are related in a manner
analogous to those arising in Theorem \ref{thm:config}.  If $\bE \cong
\C$, conditions (1) and (2) imply that $X_{\ell}(\bE)$ is a $K(G,1)$
space, where $G=\pi_1(X_{\ell}(\bE))$ is the fundamental group of
$X_{\ell}(\bE)$, as is readily seen from the homotopy sequence of a
bundle.  In this case, condition (3) further implies that the group
$G$ admits the structure of an iterated semidirect product of free
groups, and condition (4) restricts the type of free group
automorphisms arising in this structure.  These conditions determine
the additive structure of the Lie algebra $E^*_0(G)$, see \cite{FR1}
and Section \ref{sec:DCS}.  For higher dimensional $\bE$, conditions
(1)--(3) imply that the homology of the loop space of $X_{\ell}(\bE)$
is isomorphic to the universal enveloping algebra of the Lie algebra
$\pi_*(\Omega X_{\ell}(\bE))/\Tors$, and determine the additive
structure of $\Prim H_*(\Omega X_{\ell}(\bE);\Z)$, see \cite{CX} and
Section \ref{sec:loop}.  For higher dimensional $\bE$, these
conditions have analogous implications for the homology of an iterated
loop space $\Omega^q X_\ell(\bE)$ with $q>1$, and the Poisson algebra
structure admitted by this homology, see \cite{CLM} and
Section~\ref{sec:poisson}.

A brief indication how one may analyze and compare the Lie algebras
arising for various choices of $\bE$ is given next.  First, there is a
variant of the classical Freudenthal suspension, relating reduced
suspensions and loop spaces as indicated below, where the maps are
induced by (homology) suspensions.
\begin{equation*} \label{eq:indecomposables}
\begin{CD}
H_{2k-2}(\Omega X_\ell(\C^k)) @<<< H_{2k-1}(X_\ell(\C^k)) \\
@.                            @VVV \\
@.                            H_{2k}(\Sigma X_\ell(\C^k)) @>>>
H_{2k}(\Omega X_\ell(\C^{k+1}))
\end{CD}
\end{equation*}
If $k\ge 2$, conditions (1)--(3) above imply that these maps are all
(additive) isomorphisms.  In the case $k=1$, these maps yield an
additive isomorphism $E^1_0(G) = H_1(X_\ell(\C)) \cong H_{2}(\Omega
X_\ell(\C^2))$ where $G=\pi_1(X_\ell(\C))$.  While this comparison
does not in general preserve the structures of these Lie algebras, it
does provide a geometric way to compare indecomposable elements in
these Lie algebras.

To determine the Lie algebra structure, let $S$ be a sphere of
appropriate dimension and $A:S \to X_\ell(\bE)$ a map representing
a (reduced) homology generator in minimal degree.  Consider the
pullback $\xi(\bE)$ of the bundle $X_{\ell+1}(\bE)\to X_\ell(\bE)$
along the map $A$:
\begin{equation*} \label{eq:pullback}
\begin{CD}
\bE\setminus\cQ_\ell(\bE) @>>> \xi(\bE)        @>>> S\\
@|                             @VVV                 @VV{A}V\\
\bE\setminus\cQ_\ell(\bE) @>>> X_{\ell+1}(\bE) @>>> X_\ell(\bE)
\end{CD}
\end{equation*}
These bundles admit compatible cross-sections by condition (3).  
There is consequently a morphism of extensions of Lie algebras
\begin{equation*} \label{eq:morphism}
\begin{CD}
0 @>>> \cL(\bE\setminus\cQ_\ell(\bE)) @>>> \cL(\xi(\bE)) @>>> \cL(S)
@>>>0\\
@.     @VV{\id}V                         @VVV
@VV{A_*}V\\
0 @>>> \cL(\bE\setminus\cQ_\ell(\bE)) @>>> \cL(X_{\ell+1}(\bE)) @>>>
\cL(X_\ell(\bE)) @>>>0
\end{CD}
\end{equation*}
where $\cL(\bullet)$ denotes the Lie algebra obtained from the
descending central series of the fundamental group if $\bE\cong \C$,
and the graded Lie algebra of primitive elements in the homology of
the loop space for higher dimensional $\bE$.  Knowledge of the
extension $0\to \cL(\bE\setminus\cQ_\ell(\bE)) \to \cL(\xi(\bE)) \to
\cL(S) \to 0$ and the map $A_*:\cL(S) \to \cL(X_\ell(\bE))$ for all
homology generators completely determines the structure of the Lie
algebra $\cL(X_{\ell+1}(\bE))$.  In favorable situations, one can show
that the extensions of Lie algebras which arise as $\bE$ varies are,
apart from grading, isomorphic by carefully combining these
considerations with the aforementioned comparison of indecomposables.

Several natural families of examples which fit in the framework
described above are given next.  These examples either may be or have
been studied using (variants of) the techniques sketched above.  Let
$M$ be a manifold, and $\Gamma$ a group which acts properly
discontinuously on $M$.  The orbit configuration space
$F_{\Gamma}(M,\ell)$ consists of all $\ell$-tuples of points in $M$,
no two of which lie in the same $\Gamma$-orbit.

First, consider orbit configuration spaces of the form
$F_{\Gamma}(\bE\times \C^{n},\ell)$, where $\Gamma$ operates
diagonally of $\bE\times\C^{n}$, and trivially on $\C^{n}$.  
Relevant examples include the following.
\begin{enumerate}
    \item[(a)] A parameterized lattice $\Gamma$ acting on $\bE=\C$, so
    that the orbit space is an elliptic curve.  The orbit
    configuration spaces associated to the action of the standard
    integral lattice were the subject of \cite{CX}, where it is shown
    that the analogue of Theorem \ref{thm:config} holds for these
    spaces.

\smallskip

    \item[(b)] A discrete group $\Gamma$ acting properly
    discontinuously on the upper half-plane $\bE=\bH$, so that the
    orbit space is a complex curve.

\smallskip

    \item[(c)] A torsion free subgroup of $\Gamma<\Symp(2g,\Z)$ acting
    properly discontinuously on Siegel upper half-space $\bE=\bH^{g}$.

\smallskip

    \item[(d)] A torsion free subgroup $\Gamma$ of the mapping class
    group for genus $g$ surfaces, acting on Teichmuller space $\bE$.
\end{enumerate}

Second, let $M=\C^{k}\setminus\{0\}$ and let $\Gamma=\Z/p\Z$ act
freely on $M$ by rotations.  The orbit configuration spaces
$F_{\Gamma}(M,\ell)$ were the subject of \cite{Xi} and \cite{DCmono},
the results of which combine to show that the analogue of Theorem
\ref{thm:config} also holds for these spaces.

In the instances where conditions (1)--(4) hold, one obtains
generalizations of the universal Yang-Baxter Lie algebra,
parameterized by the group $\Gamma$.  This is the case for the family
$F_{\Gamma}(M,\ell)$ of orbit configuration spaces above, where
$M=\C^k\setminus\{0\}$ and $\Gamma=\Z/p\Z$.  As noted by D.~Matei, the
resulting generalized Yang-Baxter Lie algebra with cyclic symmetry is
of use in constructing Vassiliev invariants of links in the lens space
$L(p,1)$.  The Lie algebras arising from other families of orbit
configuration spaces may be of similar use for other three-manifolds,
among other potential applications.  The orbit configuration spaces
$F_{\Z/p\Z}(\C^k\setminus\{0\},\ell)$, and the classical configuration
spaces $F(\C^k,\ell)$, may be realized as complements of finite
hyperplane or subspace arrangements.  This led to speculation in
\cite{CX} that similar results may hold for fiber-type arrangements
whose complements, like configuration spaces, admit iterated bundle
structure.

Let $\A$ be a hyperplane arrangement in $\C^{\ell}$, a finite
collection of codimension one affine subspaces, with complement
$M(\A)=\C^{\ell}\setminus \bigcup_{H\in\A} H$.  See Orlik and Terao
\cite{OT} as a general reference on arrangements.  Given a hyperplane
$H\subset\C^\ell$, let $H^{k}$ be the codimension $k$ affine subspace
of $\C^{k\ell}=(\C^{\ell})^{k}$ consisting of all $k$-tuples of points
in $\C^{\ell}$, each of which lies in $H$.  For each positive integer
$k$, the elements of the hyperplane arrangement $\A$ may be used in
this way to obtain an arrangement $\A^k$ of complex codimension $k$
subspaces in $\C^{k\ell}$, with complement
$M(\A^k)=\C^{k\ell}\setminus \bigcup_{H\in\A} H^k$.

When $\A$ is a fiber-type hyperplane arrangement, the behavior of the
family of spaces $\{X_{\ell}(\C^{k})=M(\A^k),k\ge 1\}$ is reminiscent
of that of the family $\{F(\C^k,n),k\ge 1\}$ of configuration spaces. 
Let $G=\pi_1(M(\A))$ be the fundamental group of the complement of the
fiber-type arrangement $\A$ in $\C^\ell$, and let $E^{*}_{0}(G)$ be
the Lie algebra obtained from the descending central series of $G$. 
The main result of this article is as follows.

\begin{thm} \label{thm:CXconj}
For $k\ge 1$, the homology of $\Omega M(\A^{k+1})$, the loop space of
the complement of the subspace arrangement $\A^{k+1}$ in
$\C^{(k+1)\ell}$, is isomorphic to the universal enveloping algebra of
the graded Lie algebra $E^*_0(G)_k$.  Moreover,
\begin{enumerate}
\item[(a)] The image of the Hurewicz homomorphism
\[
\pi_*(\Omega M(\A^{k+1})) \to H_*(\Omega M(\A^{k+1});\Z)
\]
is isomorphic to $E^*_0(G)_k$; and

\smallskip

\item[(b)] The Hurewicz homomorphism induces isomorphisms of graded
Lie algebras
\[
\pi_*(\Omega M(\A^{k+1}))/\Tors \to \Prim H_*(\Omega
M(\A^{k+1});\Z) \cong E^*_0(G)_k,
\]
where the Lie algebra structure of the source is induced by the
Samelson product.

\end{enumerate}
\end{thm}

The paper is organized as follows.
\begin{enumerate}
\item[{\S}\ref{sec:arrangements}.]  Given a hyperplane arrangement
$\A\subset\C^\ell$, there is an associated arrangement of codimension
$k$ subspaces $\A^k\subset\C^{k\ell}$.  The combinatorics and topology
of the subspace arrangement $\A^k$ are studied in this section.

\smallskip

\item[{\S}\ref{sec:fiber-type}.]  The topology of the subspace
arrangement $\A^k$, in the instance where the underlying hyperplane
arrangement $\A$ is fiber-type, is further studied in this section.

\smallskip

\item[{\S}\ref{sec:DCS}.]  The (known) structure of the Lie algebra
$E^*_0(G)$ associated to the descending central series of the
fundamental group $G=\pi_1(M(\A))$ of the complement of a fiber-type
hyperplane arrangement $\A$ is analyzed in this section.

\smallskip

\item[{\S}\ref{sec:loop}.]  The structure of the Lie algebra of
primitive elements in the homology of the loop space of the complement
of the subspace arrangement $\A^k$ is analyzed in this section, and
the isomorphisms of graded Lie algebras asserted in
Theorem~\ref{thm:CXconj} are established. 

\smallskip

\item[{\S}\ref{sec:poisson}.]  The Poisson algebra structure on the
homology of an iterated loop space of the complement of the subspace
arrangement $\A^k$ is briefly analyzed in this section.
\end{enumerate}

\section{Redundant Arrangements} \label{sec:arrangements}

Let $H$ be an affine hyperplane in $\C^\ell$, an affine subspace of
codimension one.  For each positive integer $k$, there is an affine
subspace $H^k$ of codimension $k$ in $\C^{k\ell}$ obtained from $H$ in
the following manner.  Choose coordinates $\b{x}=(x_1,\dots,x_\ell)$
on $\C^\ell$, and $(\b{x}_1,\dots,\b{x}_k)$ on
$\C^{k\ell}=\C^\ell\times\dots\times\C^\ell$, where for each $i$,
$\b{x}_i=(x_{i,1},\dots,x_{i,\ell})\in\C^\ell$.  Then, if the
hyperplane $H$ in $\C^\ell$ is given by $H=\bigl\{\b{x}\in\C^\ell\mid
\sum_{j=1}^{\ell}a_j x_j = b\bigr\}$, define a codimension $k$ affine
subspace $H^k$ in $\C^{k\ell}$ by
$H^k=\bigl\{(\b{x}_1,\dots,\b{x}_k)\in\C^{k\ell}\mid
\sum_{j=1}^{\ell}a_j x_{i,j} = b,\ 1\le i\le k\bigr\}$.

Now let $\A$ be a hyperplane arrangement in $\C^\ell$, a finite
collection of (affine) hyperplanes.  Via the above process, there is
an arrangement $\A^k=\{H^k \mid H \in \A\}$ of codimension $k$ affine
subspaces in $\C^{k\ell}$ obtained from $\A$.  For evident reasons,
call the subspace arrangement $\A^k$ {\em redundant}.  A brief
description of the relationship between the combinatorics and topology
of the hyperplane arrangement $\A=\A^1$ and the redundant subspace
arrangement $\A^k$ is given in this section.

\begin{exm} \label{exm:braidarrangement}
Let $\A_{n}$ be the braid arrangement in $\C^{n}$, consisting of the 
hyperplanes $H_{i,j}=\{\b{x}\in\C^{n}\mid x_{i}=x_{j}\}$.  As is well 
known, the complement $M(\A_{n})=F(\C,n)$ is the configuration space 
of $n$ points in $\C$.

For each positive integer $k$, the associated redundant arrangement
$\A_{n}^{k}$ consists of subspaces
$H_{i,j}^{k}=\{(\b{x}_{1},\dots,\b{x}_{k})\in (\C^{n})^{k}\mid
x_{p,i}=x_{p,j},\ 1\le p\le k\}$.  These subspaces may be
realized as $H_{i,j}^{k}=\{(\b{y}_{1},\dots,\b{y}_{n})\in
(\C^{k})^{n}\mid \b{y}_{i}=\b{y}_{j}\}$.  Thus the complement
$M(\A^{k}_{n})=F(\C^{k},n)$ is the configuration space of $n$ points
in $\C^{k}$.
\end{exm}

For an arbitrary hyperplane arrangement $\A$, and for each $k$, let
$\sL(\A^k)$ be the intersection poset of the arrangement $\A^k$, the
partially ordered set of non-empty multi-intersections of elements of
$\A^k$.  Order $\sL(\A^k)$ by reverse inclusion, and include the
ambient space $\C^{k\ell}$ in $\sL(\A^k)$ as the minimal element,
corresponding to the intersection of no elements of $\A^k$.  For the
hyperplane arrangement $\A=\A^1$, it is known that $\sL(\A)$ is a
geometric poset, see \cite[Section 2.3]{OT}.  This need not be the
case for an arbitrary subspace arrangement.  However, for redundant
arrangements, the following holds.

\begin{prop} \label{prop:samelattice}
If $\A$ is a hyperplane arrangement, then $\sL(\A) \cong \sL(\A^k)$
for all $k$.
\end{prop}
\begin{proof}
It will be shown that the bijection between $\A$ and $\A^k$ given by
$H \leftrightarrow H^k$ induces an isomorphism of posets $\sL(\A)
\cong \sL(\A^k)$.  To establish this, it suffices to show that to each
codimension $r$ flat $X \in \sL(\A)$ there corresponds a codimension
$kr$ flat $X^k \in \sL(\A^k)$.

Write $\A=\{H_1,\dots,H_n\}$, where $H_i=\bigl\{\b{x}\in\C^\ell \mid
\sum_{j=1}^\ell a_{i,j} x_j=b_i\bigr\}$, and let $X=H_1 \cap \dots
\cap H_m$.  The flat $X$ may be realized as the set of solutions of
the system of equations $A\b{x}=\b{b}$, where $A=(a_{i,j})$ is
$m\times \ell$ and $\b{b}=(b_1,\dots,b_m)$.  Then, $X$ has codimension
$r$ in $\C^\ell$ if and only if $\rank [A \mid \b{b}] = \rank A =
\ell-r$ if and only if
\[
\rank
\left[\begin{matrix}
A& &      & & | &\b{b}\\
 &A&      & & | &\b{b}\\
 & &      & & | &   \\
 & &      &A& | &\b{b}
\end{matrix}\right] =
\rank
\left[\begin{matrix}
A& &      & \\
 &A&      & \\
 & &      & \\
 & &      &A
\end{matrix}\right] =
k(\ell-r)
\]
if and only if $X^k=H_1^k \cap \dots \cap H_m^k$ has codimension $kr$
in $\C^{k\ell}$.
\end{proof}

For each $k$, let $M(\A^{k}) = \C^{k\ell} \setminus
\bigcup_{H^{k}\in\A^{k}} H^{k}$ denote the complement of the
(subspace) arrangement $\A^{k}$.  In the case $k=1$, the cohomology of
the hyperplane complement $M(\A)=M(\A^{1})$ is well known.  
It is isomorphic to the Orlik-Solomon algebra $\sA(\A)$, see
\cite[Sections 3.1,~3.2]{OT}.  A family of algebras which includes
$\sA(\A)$ is defined next.

For each positive integer $k$, let $\E_{2k-1}[k] = \bigoplus_{H\in\A}
\Z e_{H}^k$ be a free $\Z$-module generated by degree $2k-1$ elements
$e_H^k$ in one-to-one correspondence with the hyperplanes of $\A$. 
Let $\E[k]=\bigwedge \E_{2k-1}[k]$ be the exterior algebra of
$\E_{2k-1}[k]$, and denote by $\I[k]$ the ideal of $\E[k]$ generated
by the homogeneous elements
\begin{align*}
\sum_{p=1}^q (-1)^{p-1} e^k_{H_1}\wedge\cdots
\widehat{e^k_{H_p}}\cdots
\wedge e^k_{H_q}
&\quad\text{if}\quad 0 \le \codim H_1 \cap \dots \cap H_q < q,\\
e^k_{H_1}\wedge\dots \wedge e^k_{H_q}
&\quad\text{if}\quad  H_1 \cap \dots \cap H_q = \emptyset.
\end{align*}
Let $\sA[k]=\E[k]/\I[k]$.  The Orlik-Solomon algebra is then given by
$\sA(\A)=\sA[1]$.

Proposition \ref{prop:samelattice} may be used to determine the
cohomology of $M(\A^{k})$ for $k>1$ in terms of that of $M(\A)$.  Let
$P(\A^{k},t)=\sum_{q\ge 0} b_{q}(M(\A^{k}))\cdot t^{q}$ be the
Poincar\'e polynomial of $M(\A^{k})$, where $b_{q}(X)$ is the $q$-th
Betti number of $X$.  Results of Goresky and MacPherson \cite{GM}, and
Yuzvinsky \cite{Yuz}, see also Feichtner and Ziegler \cite{FZ},
together with Proposition \ref{prop:samelattice}, yield the following.

\begin{cor} \label{cor:cohomology}
Let $\A$ be a hyperplane arrangement in $\C^{\ell}$.
\begin{enumerate}
\item For each $k$, the integral (co)homology of $M(\A^{k})$ is
torsion free, and we have $P(\A^{k},t)=P(\A,t^{2k-1})$.

\smallskip

\item For each $k$, the cohomology algebra of $M(\A^{k})$ isomorphic
to the algebra $\sA[k]$, $H^{*}(M(\A^{k});\Z) \cong \sA[k]$.
\end{enumerate}
\end{cor}

An explicit basis for the first non-zero (reduced) homology group,
$H_{2k-1}(M(\A^k);\Z)$, of the complement of the subspace arrangement
$\A^k$ is recorded next.  Let $L\subset\C^{\ell}$ be a complex line
that is transverse to the hyperplane arrangement $\A$.  Write
$L=\{t\cdot \b{u}+\b{v}\}$ where $\b{u},\b{v} \in\C^\ell$ are fixed
and $t\in\C$ varies.  For each hyperplane $H$ of $\A$, the
intersection $L\cap H$ is a point, say $\b{q}_H = \tau_H \cdot
\b{u}+\b{v}$ for some $\tau_H\in\C$.  The following is immediate.

\begin{lem}
The subspace $L^k = \{(t_1\cdot \b{u}+\b{v},\dots,t_k\cdot
\b{u}+\b{v}) \mid t_1,\dots,t_k \in \C\}$ of $\C^{k\ell}$ is
transverse to the subspace arrangement $\A^k \subset \C^{k\ell}$.  For
each subspace $H^k$ of $\A^k$, the intersection $L^k \cap H^k$ is the
point $(\b{q}_H,\dots,\b{q}_H)=(\tau_H\cdot
\b{u}+\b{v},\dots,\tau_H\cdot \b{u}+\b{v})$.
\end{lem}

Let $S^{2k-1}$ be the unit sphere in $\C^k$.  For $\epsilon > 0$
sufficiently small, the point
\[
\bigl((\tau_H+\epsilon' z_1)\cdot \b{u}+\b{v},\dots,
(\tau_H+\epsilon' z_k)\cdot \b{u}+\b{v}\bigr) \in L^k
\]
lies in the intersection $L^k \cap M(\A^k)$ for all $\epsilon'$,
$0<\epsilon' \le \epsilon$.  Fix such an $\epsilon$, and define a map
$c^k_H:S^{2k-1} \to L^k \cap M(\A^k)$ using the above formula:
\begin{equation} \label{eq:cycles}
c^k_H(\b{z})=c^k_H(z_1,\dots,z_k) =
\bigl((\tau_H+\epsilon z_1)\cdot \b{u}+\b{v},\dots,
(\tau_H+\epsilon z_k)\cdot \b{u}+\b{v}\bigr).
\end{equation}
Let $\iota_{2k-1}$ be the fundamental class of
$H_{2k-1}(S^{2k-1};\Z)$, and denote the image of
$(c^k_H)_*(\iota_{2k-1}) \in H_{2k-1}(L^k \cap M(\A^k);\Z)$ under the
map induced by the natural inclusion $L^k \cap M(\A^k) \hookrightarrow
M(\A^k)$ by $C^k_H \in H_{2k-1}(M(\A^k);\Z)$.

\begin{prop} \label{prop:basis}
The classes $\{C^k_H \mid H \in \A\}$ form a basis for
$H_{2k-1}(M(\A^k);\Z)$.
\end{prop}
\begin{proof}
For $H\in\A$, define $p^k_H:L^k\cap M(\A^k) \to S^{2k-1}$ by
\[
p^k_H(t_1\cdot \b{u}+\b{v},\dots,t_k\cdot \b{u}+\b{v}) =
\frac{\b{t} - \tau_H\cdot\b{e}}{\| \b{t} - \tau_H\cdot\b{e} \|},
\]
where $\b{t}=(t_1,\dots,t_k)$ and $\b{e}=(1,\dots,1)$ are in $\C^k$. 
It is then readily checked that $p^k_H \circ c^k_H = \id:S^{2k-1} \to
S^{2k-1}$ is the identity map.  Furthermore, if $H' \neq H$ is another
hyperplane of $\A$, the composition $p^k_H \circ c^k_{H'}$ is given by
\[
p^k_H \circ c^k_{H'}(\b{z}) =
\frac{\b{z} +\frac{1}{\epsilon}(\tau_{H'} - \tau_H)\cdot\b{e}}
{\| \b{z} +\frac{1}{\epsilon}(\tau_{H'} - \tau_H)\cdot\b{e} \|},
\]
so is null-homotopic.  Consequently, the classes
$(c^k_H)_*(\iota_{2k-1}) \in H_{2k-1}(L^k \cap M(\A^k);\Z)$ form a
basis.  Finally, using stratified Morse theory, one can show that the
relative homology group $H_{i}(M(\A^{k}),L^k \cap M(\A^k);\Z)$
vanishes for $i < 4k-2$, see \cite[Parts II,~III]{GM}.
It follows that the natural inclusion
$L^k \cap M(\A^k) \hookrightarrow M(\A^k)$ induces an isomorphism
$H_{2k-1}(L^k \cap M(\A^k);\Z) \xrightarrow{\,\sim\,}
H_{2k-1}(M(\A^k);\Z)$.  So the classes $C^k_H$ form a basis for
$H_{2k-1}(M(\A^k);\Z)$ as asserted.
\end{proof}

\begin{rem} \label{rem:dual}
The cohomology classes $(C^k_H)^* \in H^{2k-1}(M(\A^k);\Z)$ dual to
the classes $C^k_H\in H_{2k-1}(M(\A^k);\Z)$ generate the cohomology
algebra $H^*(M(\A^k);\Z)$.  Let $a_H^k \in \sA[k]$ denote the image of
$e_H^k \in \E[k]$ under the natural projection.  Then the map
$H^{2k-1}(M(\A^k);\Z) \to \sA_{2k-1}[k]$, $(C^k_H)^* \mapsto a_H^k$,
induces an isomorphism of algebras $H^*(M(\A^k);\Z)
\xrightarrow{\,\sim\,} \sA[k]$, see Corollary \ref{cor:cohomology}.
\end{rem}

\section{Linearly Fibered Arrangements} \label{sec:fiber-type}

In this section, the topology of those redundant arrangements arising
from strictly linearly fibered and fiber-type hyperplane arrangements
is studied further.  Recall the definition of arrangements of the
former type from \cite{FR1,OT}.

\begin{defn} \label{defn:slfdef}
A hyperplane arrangement $\A$ in $\C^{\ell+1}$ is {\em strictly
linearly fibered} if there is a choice of coordinates
$(\b{x},z)=(x_1,\dots,x_\ell,z)$ on $\C^{\ell+1}$ so that the
restriction, $p$, of the projection $\C^{\ell+1}\to\C^\ell$,
$(\b{x},z)\mapsto \b{x}$, to the complement $M(\A)$ is a fiber bundle
projection, with base $p(M(\A))=M(\B)$, the complement of an
arrangement $\B$ in $\C^\ell$, and fiber the complement of finitely
many points in $\C$.  Refer to the hyperplane arrangement $\A$ as
strictly linearly fibered over $\B$.
\end{defn}

The complements of hyperplane arrangements of this type are closely
related to configuration spaces, as we now illustrate.  For each
hyperplane $H$ of $\A$, let $f_H$ be a linear polynomial with $H=\ker
f_H$.  Then $Q(\A)=\prod_{H\in\A} f_H$ is a defining polynomial for
$\A$.  From the definition, if $\A$ is strictly linearly fibered over
$\B$ and $|\A|=|\B|+n$, there is a choice of coordinates for which a
defining polynomial for $\A$ factors as
\begin{equation} \label{eq:defpoly}
Q(\A) = Q(\B)\cdot \phi(\b{x},z),
\end{equation}
where $Q(\B)=Q(\B)(\b{x})$ is a defining polynomial for $\B$, and
$\phi(\b{x},z)$ is a product
\[
\phi(\b{x},z) = (z-g_1(\b{x}))(z-g_2(\b{x}))\cdots (z-g_n(\b{x})),
\]
with $g_j(\b{x})$ linear.  Define $g:\C^{\ell}\to\C^{n}$ by
\begin{equation} \label{eq:rootmap}
g(\b{x})=\bigl(g_1(\b{x}),g_2(\b{x}),\dots,g_n(\b{x})\bigr),
\end{equation}
Since $\phi(\b{x},z)$ necessarily has distinct roots for any $\b{x}\in
M(\B)$, the restriction of $g$ to $M(\B)$ takes values in the
configuration space $F(\C,n)$.  The following result was proven by the
first author, see \cite[Theorem~1.1.5, Corollary~1.1.6]{DCmono}.

\begin{thm}
\label{thm:slfpullback}
Let $\B$ be an arrangement of $m$ hyperplanes, and let $\A$ be an
arrangement of $m+n$ hyperplanes which is strictly linearly fibered
over $\B$.  Then the bundle $p:M(\A)\to M(\B)$ is equivalent to the
pullback of the bundle of configuration spaces $p_{n+1}:F(\C,n+1)\to
F(\C,n)$ along the map $g$.  Consequently, the bundle $p:M(\A)\to
M(\B)$ admits a cross-section and has trivial local coefficients in
homology.
\end{thm}
Since it is relevant to the subsequent discussion, a proof is
included.
\begin{proof}
Denote points in $F(\C,n+1)$ by $(\b{y},z)$, where
$\b{y}=(y_1,\dots,y_n)\in F(\C,n)$ and $z\in \C$ satisfies $z\neq y_j$
for each $j$.  Similarly, denote points in $M(\A)$ by $(\b{x},z)$,
where $\b{x}\in M(\B)$ and $\phi(\b{x},z)\neq 0$.  In this notation,
we have $p_{n+1}(\b{y},z)=\b{y}$ and $p(\b{x},z)=\b{x}$.  Let
$E=\bigl\{\bigl(\b{x},(\b{y},z)\bigr) \in M(\B) \times F(\C,n+1) \mid
g(\b{x})=\b{y}\bigr\}$ be the total space of the pullback of
$p_{n+1}:F(\C,n+1)\to F(\C,n)$ along the map $g$.  It is then readily
checked that the map $h:M(\A) \to E$ defined by $h(\b{x},z) =
\bigl(\b{x},(g(\b{x}),z)\bigr)$ is an equivalence of bundles.

Since the bundle $p_{n+1}:F(\C,n+1)\to F(\C,n)$ admits a
cross-section, so does the pullback $p:M(\A)\to M(\B)$.  Furthermore,
the structure group of the latter is the pure braid group $P_n$. 
Consequently, the action of the fundamental group of the base $M(\B)$
on that of the fiber $\C\setminus\{n\ \text{points}\}$ is by pure
braid automorphisms.  As such, this action is by conjugation (see for
instance \cite{Bi} or \cite{Ha}), hence is trivial in homology.
\end{proof}

It is now shown that redundant strictly linearly fibered arrangements
admit (linear) fibrations, just as their codimension one progenitors
do.

\begin{thm}
\label{thm:slfsoupup}
Let $\A$ be a hyperplane arrangement in $\C^{\ell+1}$ which is
strictly linearly fibered over $\B$, with projection $p:M(\A)\to
M(\B)$ induced by the map $\C^{\ell+1}\to\C^{\ell}$ given by
$(x_{1},\dots,x_{\ell},z)\mapsto (x_{1},\dots,x_{\ell})$.  Then for
each $k$, the map $\C^{k(\ell+1)}\to\C^{k\ell}$ given by
$(\b{x}_{1},\dots,\b{x}_{\ell},\b{z}) \mapsto
(\b{x}_{1},\dots,\b{x}_{\ell})$ induces a fiber bundle projection
$p^{k}:M(\A^{k})\to M(\B^{k})$.  Furthermore, the bundle
$p^{k}:M(\A^{k})\to M(\B^{k})$ admits a cross-section.
\end{thm}
\begin{proof}
By the previous result, the bundle $p:M(\A)\to M(\B)$ is equivalent to
the pullback of $p_{n+1}:F(\C,n+1)\to F(\C,n)$ along the map $g$ of
\eqref{eq:rootmap}.  An analogous result for the complements of the
subspace arrangements $\A^k$ and $\B^k$ is established next.  For $k
\ge 2$, view $\C^{k\ell}$ as $(\C^\ell)^k$ and $\C^{kn}$ as
$(\C^k)^n$.  Denote points in the configuration space $F(\C^k,n+1)$ by
$(\b{y}_1,\dots,\b{y}_n,\b{z})$, where $(\b{y}_1,\dots,\b{y}_n) \in
F(\C^k,n)$ and $\b{z} \neq \b{y}_j$ for each $j$.  Define
$g^k:\C^{k\ell} \to \C^{kn}$ by
\begin{equation} \label{eq:kroot}
g^k(\b{x}_1,\dots,\b{x}_k) = \Bigl(
\bigl(g_1(\b{x}_1),\dots,g_1(\b{x}_k)\bigr),\dots\,\dots,
\bigl(g_n(\b{x}_1),\dots,g_n(\b{x}_k)\bigr)\Bigr).
\end{equation}
where $\bigl(g_i(\b{x}_1),\dots,g_i(\b{x}_k)\bigr) \in \C^k$ for each
$i$.  It is readily checked that the restriction of $g^k$ to $M(\B^k)$
takes values in the configuration space $F(\C^k,n)$.  Let $\pi^k:E^k
\to M(\B^k)$ be the pullback of the bundle $p_{n+1}^k: F(\C^k,n+1) \to
F(\C^k,n)$ along this restriction, with total space $E^k$ consisting
of all points
\[
\bigl((\b{x}_1,\dots,\b{x}_k),(\b{y}_1,\dots,\b{y}_n,\b{z})\bigr) \in
M(\B^k) \times F(\C^k,n+1)
\]
for which
$g^k(\b{x}_1,\dots,\b{x}_k)=p_{n+1}^k(\b{y}_1,\dots,\b{y}_n,\b{z})=
(\b{y}_1,\dots,\b{y}_n)$.

Since the hyperplane arrangement $\A$ is strictly linearly fibered
over $\B$, the complement of the subspace arrangement $\A^k$ may be
realized as
\[
M(\A^k)=\{(\b{x}_1,\dots,\b{x}_k,\b{z}) \in M(\B^k) \times \C^k \mid
\b{z} \neq \bigl(g_i(\b{x}_1),\dots,g_i(\b{x}_k)\bigr) \ \text{for}\
1\le i\le n\}.
\]
Define $h^k:M(\A^k) \to E^k$ by $h^k(\b{x}_1,\dots,\b{x}_k,\b{z})=
\bigl((\b{x}_1,\dots,\b{x}_k),(g^k(\b{x}_1,\dots,\b{x}_k),\b{z})\bigr)$.

The map $h^k$ is a homeomorphism.
Moreover, the following diagram commutes.
\begin{equation} \label{eq:kpullback}
\begin{CD}
M(\A^k)   @>{h^k}>> E^k\\
@VV{p^k}V           @VV{\pi^k}V\\
M(\B^k)   @>{\id}>> M(\B^k)
\end{CD}
\end{equation}
It follows that $p^k:M(\A^k)\to M(\B^k)$ is a bundle which is
equivalent to the pullback of the bundle of configuration spaces
$p_{n+1}^k: F(\C^k,n+1) \to F(\C^k,n)$ along the map $g^k:M(\B^k)\to
F(\C^k,n)$, and therefore has a cross-section.
\end{proof}

An analysis of map in homology induced by the map $g^k:M(\B^k) \to
F(\C^k,n)$ defined in \eqref{eq:kroot} is given next.  For $1\le i<j
\le n$, define $p_{i,j}:F(\C^k,n) \to S^{2k-1}$ by
$p_{i,j}(\b{y}_1,\dots,\b{y}_n)=(\b{y}_j-\b{y}_i)/\|\b{y}_j-\b{y}_i\|$.
 Recall that $\iota_{2k-1} \in H_{2k-1}(S^{2k-1};\Z)$ denotes the
fundamental class.  The classes $p_{i,j}^*(\iota_{2k-1})$ form a basis
for $H^{2k-1}(F(\C^k,n))$, and generate the cohomology algebra
$H^*(F(\C^k,n))$, see \cite{CLM,CT}.  Denote the dual classes in
$H_{2k-1}(F(\C^k,n)$ by $A_{i,j}$, $1\le i<j \le n$.  Note that the
classes $A_{i,j}$ may be represented by spheres linking the subspaces
$H_{i,j}^k = \{\b{y}_i = \b{y}_j\}$ in $\C^{kn}$, as in
\eqref{eq:cycles}.

As in Section \ref{sec:arrangements}, let $L=\{t\cdot\b{u}+\b{v}\}
\subset \C^\ell$ be a line transverse to the hyperplane arrangement
$\B$, and $L^k$ the corresponding codimension $k$ subspace of
$\C^{k\ell}$, transverse to the subspace arrangement $\B^k$.  Recall
the maps $c^k_H:S^{2k-1}\to L^k\cap M(\B^k)$ from \eqref{eq:cycles},
and the resulting basis $\{C^k_H \mid H \in \B\}$ for
$H_{2k-1}(M(\B^k))$ exhibited in Proposition~\ref{prop:basis}.

\begin{prop} \label{prop:inducedmap}
Let $\B \subset \C^{\ell}$ be an arrangement of complex hyperplanes,
and let $g:\C^{\ell}\to\C^{n}$ be an affine transformation whose
restriction, $g:M(\B)\to F(\C,n)$, to the complement of $\B$ takes
values in the configuration space $F(\C,n)$.  Then for every $k\ge 1$,
the induced map $(g^k)_*:H_{2k-1}(M(\B^{k});\Z) \to
H_{2k-1}(F(\C^{k},n);\Z)$ is given by $(g^k)_*(C^k_H) = \sum A_{i,j}$
for each hyperplane $H$ of $\B$, where the sum is over all distinct
$i$ and $j$ for which $g(H)$ is contained in the hyperplane
$H_{i,j}=\{y_i=y_j\}$ in $\C^n$.
\end{prop}
\begin{proof}
For each hyperplane $H$ of $\B$, let $\tilde{c}^k_H:S^{2k-1} \to
M(\B^k)$ denote the composition of $c^k_H:S^{2k-1} \to L^k \cap
M(\B^k)$ and the natural inclusion $L^k \cap M(\B^k) \hookrightarrow
M(\B^k)$.  It will be shown that the composition $p_{i,j} \circ g^k
\circ \tilde{c}^k_H: S^{2k-1} \to S^{2k-1}$ induces the identity in
homology if $g(H) \subset H_{i,j}$, and induces the trivial
homomorphism if $g(H) \not\subset H_{i,j}$, thereby establishing the
result.

For $\b{x}\in\C^{\ell}$, write $g(\b{x}) =
(g_1(\b{x}),\dots,g_n(\b{x}))$ as in \eqref{eq:rootmap}.  Then
$g^k:\C^{k\ell}\to\C^{kn}$ is given by $g^k(\b{x}_1,\dots,\b{x}_k) =
(\b{y}_1,\dots,\b{y}_n)$, where
$\b{y}_i=(g_i(\b{x}_1),\dots,g_i(\b{x}_k))$, see \eqref{eq:kroot}. 
Since the restriction of $g$ to $M(\B)$ takes values in $F(\C,n)$, the
restriction of $g^k$ to $M(\B^k)$ takes values in $F(\C^k,n)$.

 From \eqref{eq:cycles}, the map $\tilde{c}^k_H:S^{2k-1}\to M(\B^k)$
is given by $\tilde{c}^k_H(\b{z})=(\b{w}_1,\dots,\b{w}_k)$, where
$\b{w}_j=(\tau_H+\epsilon z_j)\cdot \b{u}+\b{v}$, and $L\cap H$ is the
point $\b{q}_H=\tau_H\cdot\b{u}+\b{v}$.  Let $\a_i = g_i(\b{q}_H)$,
and define $\beta_i$ by the equation
\[
g_i(\b{w}_j)=g_i((\tau_H+\epsilon z_j)\cdot \b{u}+\b{v})=
g_i(\b{q}_H + \epsilon z_j \cdot\b{u})= \a_i+\epsilon \beta_i z_j.
\]
Then, a calculation yields $g^k \circ \tilde{c}^k_H(\b{z}) =
(\a_1\cdot \b{e} +\epsilon \beta_1\cdot \b{z},\dots, \a_n\cdot \b{e}
+\epsilon \beta_n\cdot \b{z})$ and
\[
p_{i,j} \circ g^k \circ \tilde{c}^k_H(\b{z}) =
\frac{\epsilon(\beta_j-\beta_i) \b{z} + (\a_j-\a_i)\b{e}}
{\|\epsilon(\beta_j-\beta_i) \b{z} + (\a_j-\a_i)\b{e}\|},
\]
where, as before, $\b{e}=(1,\dots,1)$.  Recall that $\epsilon>0$ was
chosen sufficiently small so as to insure that the point
$(\b{w}'_1,\dots,\b{w}'_k)$, where $\b{w}'_j=(\tau_H+\epsilon'
z_j)\cdot \b{u}+\b{v}$, lies in $L^k\cap M(\B^k)$ for all $\epsilon'$,
$0<\epsilon'\le\epsilon$.  Since $g^k:M(\B^k)\to F(\C^k,n)$, it
follows that $g^k \circ \tilde{c}^k_H(\b{z}) \in F(\C^k,n)$ for all
$\b{z} \in S^{2k-1}$.  In other words, $\epsilon(\beta_j-\beta_i)
\b{z} + (\a_j-\a_i)\b{e} \neq \b{0}$ for all distinct $i$ and $j$.  

If $g(H) \not\subset H_{i,j}$, then $g(\b{q}_H) \notin H_{i,j}$ since
$\b{q}_H = L\cap H$ is generic in $H$.  Thus, $\a_i = g_i(\b{q}_H)
\neq g_j(\b{q}_H) = \a_j$, and the point $(\a_i\b{e}+\epsilon'
\beta_i\b{z},\a_j\b{e}+\epsilon' \beta_j\b{z})$ lies in the
configuration space $F(\C^k,2)$ for all $\epsilon'\le \epsilon$,
including $\epsilon'=0$.  It follows that $p_{i,j} \circ g^k \circ
c^k_H$ is trivial in homology in this instance.

If, on the other hand, $g(H) \subset H_{i,j}$, then $\a_i = \a_j$ and
thus $\beta_j-\beta_i$ is necessarily non-zero.  In this instance,
$p_{i,j} \circ g^k \circ \tilde{c}^k_H(\b{z}) = \lambda\cdot \b{z}$,
where $\lambda \in S^1 \subset\C^*$ is given by
$\lambda=(\beta_j-\beta_i)/|\beta_j-\beta_i|$, which clearly induces
the identity in homology.
\end{proof}

These results extend immediately to fiber-type arrangements, defined
next.

\begin{defn} \label{defn:ftdef}
An arrangement $\A=\A_1$ of finitely many points in $\C^1$ is {\em
fiber-type}.  An arrangement $\A=\A_\ell$ of hyperplanes in $\C^\ell$
is {\em fiber-type} if $\A$ is strictly linearly fibered over a
fiber-type hyperplane arrangement $\A_{\ell-1}$ in $\C^{\ell-1}$.
\end{defn}

Let $\A$ be a fiber-type hyperplane arrangement in $\C^\ell$.  Then
there is a choice of coordinates $(x_{1},\dots,x_{\ell})$ on
$\C^{\ell}$ so that a defining polynomial for $\A$ factors as
$Q(\A)=\prod_{j=1}^{\ell} Q_{j}(x_{1},\dots,x_{j})$, see
\eqref{eq:defpoly}.  Write
$Q_j=\prod_{i=1}^{d_j}\bigl(x_j-g_{i,j}(x_1,\dots,x_{j-1})\bigr)$,
where $d_j$ is the degree of $Q_j$ and each $g_{i,j}$ is linear.  The
non-negative integers $\{d_1,\dots,d_\ell\}$ are called the exponents
of $\A$.  For each $j\le \ell$, the polynomial $\prod_{i=1}^{j} Q_{i}$
defines a fiber-type arrangement $\A_{j}$ in $\C^{j}$ with exponents
$\{d_1,\dots,d_j\}$, and $\A_{j}$ is strictly linearly fibered over
$\A_{j-1}$.  Furthermore, the map
$g_j=(g_{1,j},\dots,g_{d_j,j}):\C^{j-1} \to \C^{d_j}$ gives rise to
maps $g_j^k:M(\A^k_{j-1}) \to F(\C^k,d_j)$ for each $k$.  Theorems
\ref{thm:slfpullback} and \ref{thm:slfsoupup} yield

\begin{thm} \label{thm:ftbundle}
Let $\A$ be a fiber-type hyperplane arrangement in $\C^\ell$ with
defining polynomial $Q(\A)=\prod_{j=1}^\ell Q_j$.  Then, for each $j$,
$2\le j \le \ell$, and each $k\ge 1$, the projection
$\C^j\to\C^{j-1}$, $(x_1,\dots,x_{j-1},x_j) \mapsto
(x_1,\dots,x_{j-1})$, gives rise to a bundle map $p_j^k:M(\A^k_{j})
\to M(\A^k_{j-1})$.  This bundle is equivalent to the pullback of the
bundle of configuration spaces $F(\C^k,d_j+1) \to F(\C^k,d_j)$ along
the map $g_j^k:M(\A^k_{j-1}) \to F(\C^k,d_j)$.  Consequently, the
bundle $p_j^k:M(\A^k_{j}) \to M(\A^k_{j-1})$ admits a cross-section,
has trivial local coefficients in homology, and the fiber is the
complement of $d_j$ points in $\C^k$.
\end{thm}

Proposition \ref{prop:inducedmap} also extends to fiber-type
arrangements.  The specific statement is omitted.

\section{The Descending Central Series} \label{sec:DCS}

In this section, the structure of the Lie algebra $E^*_0(G)$
associated to the descending central series of the fundamental group
$G$ of the complement of a fiber-type arrangement is analyzed. 
Additively, this structure is given by well known results of Falk and
Randell \cite{FR1,FR2} stated below.  Moreover, as shown by Jambu and
Papadima \cite{JP}, this Lie algebra is isomorphic to the (integral)
holonomy Lie algebra of the arrangement $\A$ defined by Kohno
\cite{K1}.  An alternate description of $E^*_0(G)$, which will
facilitate comparison with the Lie algebra of primitives in the
homology of the loop space of the complement of the subspace
arrangement $\A^k$ in Section \ref{sec:loop}, is given here.

\begin{exm} \label{exm:Pn}
Let $P_n$ be the Artin pure braid group, the fundamental group of the
configuration space $F(\C,n)$.  The structure of the Lie algebra
$E^*_0(P_n)$ was first determined rationally by Kohno \cite{Ko}.  As
observed by many authors, the following description holds over the
integers as well.  
For each $n \ge 2$, let $L[n]$ be the free Lie algebra generated by
elements $A_{1,n+1},\dots,A_{n,n+1}$.  Then the Lie algebra
$E^*_0(P_n)$ is additively isomorphic to the direct sum
$\bigoplus_{j=1}^{n-1} L[j]$, and the Lie bracket relations in
$E^*_0(P_n)$ are the infinitesimal pure braid relations, given by
\begin{equation} \label{eq:braidrels}
\begin{aligned}
{[}A_{i,j}+A_{i,k}+A_{j,k},A_{m,k}]&=0 \quad
\text{for $m=i$ or $m=j$, and}\\
[A_{i,j},A_{k,l}]&=0
\quad\text{for $\{i,j\}\cap\{k,l\}=\emptyset$.}
\end{aligned}
\end{equation}

Note that this description realizes the Lie algebra $E^*_0(P_{n+1})$
as the semidirect product of $E^*_0(P_{n})$ by $L[n]$ determined by
the Lie homomorphism $\theta_n:E^*_0(P_n) \to \Der(L[n])$ given by
$\theta_n(A_{i,j}) = \ad(A_{i,j})$.  From the
infinitesimal pure braid relations, one has
\begin{equation} \label{eq:adAij}
\ad(A_{i,j})(A_{m,n+1})=
\begin{cases}
[A_{m,n+1},A_{i,n+1}+A_{j,n+1}]&\text{if $m=i$ or $m=j$,}\\
0&\text{otherwise.}
\end{cases}
\end{equation}
This extension of Lie algebras arises topologically from the bundle
of configuration spaces $F(\C,n+1)\to F(\C,n)$.
\end{exm}

The additive structure noted above may be obtained from the
following result of Falk and Randell \cite{FR1, FR2}.

\begin{thm} \label{thm:FalkRandell}
Let $1 \to H \to G \to K \to 1$ be a split extension of groups such
that the conjugation action of $K$ is trivial on $H_1/H_2$.  Then
there is a short exact sequence of Lie algebras $0 \to E^*_0(H) \to
E^*_0(G) \to E^*_0(K) \to 0$ which is split as a sequence of abelian
groups.  Furthermore, if the descending central series quotients of
$H$ and $K$ are free abelian, then so are those of $G$.
\end{thm}

Now let $\A=\A_\ell$ be a fiber-type hyperplane arrangement in
$\C^\ell$.  The complement of $\A_\ell$ sits atop a tower of fiber
bundles
\begin{equation*} \label{eq:BundleTower}
M(\A_\ell) \xrightarrow{\ p_{\ell}\ } M(\A_{\ell-1})
\xrightarrow{p_{\ell-1}} \cdots
\xrightarrow{\ p_{2}\ } M(\A_1) = \C \setminus \{d_1
\text{ points}\},
\end{equation*}
where the fiber of $p_{j}$ is homeomorphic to the complement of
$d_{j}$ points in $\C$.  For simplicity, write $\B=\A_{\ell-1}$ and
$n=d_\ell$.  Then $\A$ is strictly linearly fibered over $\B$, and by
Theorem \ref{thm:slfpullback}, the bundle $p:M(\A)\to M(\B)$ is
equivalent to the pullback of the configuration space bundle
$p_{n+1}:F(\C,n+1)\to F(\C,n)$ along the map $g$ of
\eqref{eq:rootmap}.

Application of the homotopy exact sequence of a bundle (and induction)
shows that $M(\A)$ is a $K(G,1)$ space, where
$G=G(\A)=\pi_{1}(M(\A))$.  In light of Theorem \ref{thm:slfpullback},
there is also a commutative diagram
\begin{equation} \label{eq:groups}
\begin{CD}
1 @>>> \F_{n}    @>>> G(\A)   @>>> G(\B) @>>>1\\
@.     @VV{\id}V      @VVV         @VV{g_{\#}}V\\
1 @>>> \F_{n}    @>>> P_{n+1} @>>> P_{n} @>>>1
\end{CD}
\end{equation}
where $g_\#:G(\B)\to P_n$ is induced by $g:M(\B)\to F(\C,n)$, and the
fundamental group of the fiber $\C\setminus\{n\ \text{points}\}$ is
identified with the free group $\F_n$ on $n$ generators.  Since the
underlying bundles admit cross-sections, the rows in the diagram above
are split exact.

\begin{thm} \label{thm:DCSadditive}
Let $\A$ be a fiber-type hyperplane arrangement.  If the exponents of
$\A$ are $\{d_1,\dots,d_\ell\}$, then
$E^*_0(G(\A)) \cong L[d_1]\oplus \cdots \oplus L[d_\ell]$ as abelian
groups.
\end{thm}
\begin{proof}
The proof is by induction on $\ell$.

In the case $\ell=1$, $\A$ is an arrangement of $d=d_1$ points in
$\C$, the fundamental group of the complement is $\F_d$, the free
group on $d$ generators, and it is well known that $E^*_0(\F_d)$ is
isomorphic to the free Lie algebra $L[d]$, see for instance
\cite[Chapter IV]{Serre}.

In general, assume that the fiber-type arrangement $\A$ is strictly
linearly fibered over $\B$ and that $d_\ell=n$ as above.  Then
there is a split, short exact sequence of groups
$1 \to \F_n \to G(\A) \to G(\B) \to 1$, and by Theorem
\ref{thm:slfpullback}, the action of $G(\B)$ on $\F_n$ is by pure
braid automorphisms.  As such, this action is by conjugation, hence is
trivial on $H_*(\F_n;\Z)$.  By Theorem \ref{thm:FalkRandell}, the
descending central series quotients of $G(\A)$ are free abelian, and
there is a short exact sequence of Lie algebras
\begin{equation} \label{eq:extension}
0 \to E^*_0(\F_n) \to E^*_0(G(\A)) \to E^*_0(G(\B)) \to 0,
\end{equation}
which splits as a sequence of abelian groups.  The result follows by
induction.
\end{proof}

The additive decomposition provided by this result does not, in
general, preserve the underlying Lie algebra structure.  An inductive
description of the Lie algebra structure of $E_{0}^{*}(G(\A))$ is
given next.

\begin{thm} \label{thm:DCSstructure}
Let $\A$ and $\B$ be fiber-type hyperplane arrangements with
$|\A|=|\B|+n$, and suppose that $\A$ is strictly linearly fibered over
$\B$.  Then the Lie algebra $E^{*}_{0}(G(\A))$ is isomorphic to the
semidirect product of $E^{*}_{0}(G(\B))$ by $L[n]$ determined by the
Lie homomorphism $\Theta = \theta_{n} \circ g_{*}:E^{*}_{0}(G(\B))\to
\Der(L[n])$, where $g_{*}:E^{*}_{0}(G(\B))\to E^{*}_{0}(P_{n})$ is
induced by the map $g:M(\B)\to F(\C,n)$, and $\theta_n:E^*_0(P_n) \to
\Der(L[n])$ is given by $\theta_n(A_{i,j}) = \ad(A_{i,j})$.
\end{thm}
\begin{proof}
 From the exact sequence of Lie algebras \eqref{eq:extension}
noted above, it follows that $E^{*}_{0}(G(\A))$ is isomorphic to the
semidirect product of $E^{*}_{0}(G(\B))$ by $L[n]$ determined by the
Lie homomorphism $\Theta:E^{*}_{0}(G(\B))\to \Der(L[n])$ given by
$\Theta(b)=\ad_{L[n]}(b)$ for $b\in E^{*}_{0}(G(\B))$.  It suffices to
show that the homomorphism $\Theta$ factors as asserted.

 From the diagram \eqref{eq:groups}, and the results of Falk and
Randell stated in Theorem \ref{thm:FalkRandell}, there is
a commutative diagram of Lie algebras
with split exact rows
\[
\begin{CD}
0 @>>> L[n]    @>>> E^{*}_{0}(G(\A))   @>>> E^{*}_{0}(G(\B)) @>>>0\\
@.     @VV{\id}V    @VVV                    @VV{g_{*}}V\\
0 @>>> L[n]    @>>> E^{*}_{0}(P_{n+1}) @>>> E^{*}_{0}(P_{n}) @>>>0
\end{CD}
\]
Via the splittings, view $E^{*}_{0}(G(\B))$ and $E^{*}_{0}(P_{n})$ as
Lie subalgebras of $E^{*}_{0}(G(\A))$ and $E^{*}_{0}(P_{n+1})$
respectively.  Then for $a \in L[n]$ and $b \in E^{*}_{0}(G(\B))$, we
have $[b,a]=[g_*(b),a]$ in $L[n]$.  Thus $\ad_{L[n]}(b) =
\ad_{L[n]}(g_*(b))$ in $\Der(L[n])$ and $\Theta=\theta_n\circ g_*$.
\end{proof}

This result, together with Proposition \ref{prop:inducedmap}, provides
an inductive description of the Lie bracket structure of
$E^*_0(G(\A))$.  Recall the basis $\{C^1_H \mid H\in \B\}$ for
$H_1(M(\B);\Z) = E^1_0(G(\B))$ exhibited in Proposition
\ref{prop:basis}, and recall that the free Lie algebra $L[n]$ is
generated by $A_{1,n+1},\dots,A_{n,n+1}$.

\begin{cor} \label{cor:brackets}
For generators $C^1_H$ of $E^1_0(G(\B))$ and $A_{m,n+1}$ of $L[n]$,
one has
\[
\Theta(C^1_H)(A_{m,n+1}) = \sum_{g(H) \subset H_{i,j}}
[A_{i,j},A_{m,n+1}].
\]
\end{cor}
\begin{proof}
By Proposition \ref{prop:inducedmap}, one has $g_{*}(C^{1}_{H})= 
\sum A_{i,j}$, where the sum is over all $i$ and $j$ for which $g(H)
\subset H_{i,j}$.  The result follows.
\end{proof}

This corollary can be used to explicitly record the Lie bracket
relations in $E^*_0(G(\A))$, and to show that these relations are
combinatorial, that is, completely determined by the intersection
poset $\sL(\A)$.  The Lie algebra $E^*_0(G(\A))$ is generated by
$\{C^1_H \mid H\in \A\}$.  For a flat $X\in\sL(\A)$, write
$C^{1}_{X}=\sum_{X\subset H} C^{1}_H$.  The following was proven 
by Jambu and Papadima \cite{JP}, see also \cite{DCmono}.

\begin{thm} \label{thm:LieBrackets}
If $\A$ is a fiber-type hyperplane arrangement with exponents
$\{d_{1},\dots,d_{\ell}\}$, then the Lie bracket relations in
$E^*_0(G(\A))$ are given by
\[
[C^{1}_X,C^{1}_H]=0,
\]
for codimension two flats $X\in\sL(\A)$ and hyerplanes $H\in\A$ 
containing $X$.
\end{thm}
\begin{proof}  The proof is by induction on $\ell$.

In the case $\ell=1$, there is nothing to show since $G(\A)$ is a 
free group on $d=d_{1}$ generators, $E^*_0(G(\A))$ is isomorphic 
to the free Lie algebra $L[d]$, and there are no codimension two 
flats in $\sL(\A)$.

In general, assume that $\A$ is strictly linearly fibered over $\B$
and that $d_{\ell}=n$ as before.  Then $\A$ has a defining polynomial
of the form $Q(\A)=Q(\B) \cdot \prod_{j=1}^{n}(z-g_j(\b{x}))$, see
\eqref{eq:defpoly}.  View $\B$ as a subarrangement of
$\A=\{H \mid H\in\B\}\cup \{H_j \mid 1\le j\le n\}$, where
$H_j=\ker(z-g_j(\b{x}))$.  Then the set $\{C_H^1 \mid H\in\B\} \cup
\{C^1_{H_j} \mid 1\le j\le n\}$ generates $E^*_0(G(\A))$, where the
generators $C^1_{H_j}$ correspond to the hyperplanes $H_j$ of
$\A\setminus\B$, and to the generators $A_{j,n+1}$ of the free Lie
algebra $L[n]$ under the additive isomorphism $E^*_0(G(\A)) \cong
E^*_0(G(\B)) \oplus L[n]$.

By Theorem \ref{thm:DCSstructure}, $E^{*}_{0}(G(\A))$ is isomorphic to
an extension of $E^{*}_{0}(G(\B))$ by $L[n]$.  Consequently, the Lie
bracket relations in $E^{*}_{0}(G(\A))$ consist of those of
$E^{*}_{0}(G(\B))$, and those arising from the extension.  By
induction, the Lie bracket relations in $E^{*}_{0}(G(\B))$ are given
by $[C^1_X,C^1_H]=0$ for codimension two flats $X$ contained only in
hyperplanes $H\in\B\subset\A$.  So it remains to analyze those
relations in $E^*_0(G(\A))$ arising from the extension.  These are
given implicitly in Corollary~\ref{cor:brackets}.

Recall from \eqref{eq:braidrels} that
$[A_{i,j},A_{m,n+1}]=[A_{m,n+1},A_{i,n+1}+A_{j,n+1}]$ if
$m\in\{i,j\}$, and is zero otherwise.  Thus the results of Corollary
\ref{cor:brackets} may be recorded as
\[
\Theta(C^1_H)(A_{m,n+1}) = [C^1_H,A_{m,n+1}] = \sum_{g(H) \subset
H_{i,j}} [A_{m,n+1},(\delta_{i,m}+\delta_{j,m})(A_{i,n+1}+A_{j,n+1})],
\]
where $C^1_H \in E^*_0(G(\B)) \subset E^*_0(G(\A))$ and $\delta_{i,m}$
is the Kronecker delta.  Note that the expression on the right lies in
$L[n]$.  Under the above identifications, these relations take the
form
\begin{equation} \label{eq:relation}
[C^1_H,C^1_{H_m}] = \sum_{g(H) \subset H_{i,j}}
[C^1_{H_m},(\delta_{i,m}+\delta_{j,m})(C^1_{H_i}+C^1_{H_j})]
\end{equation}
Now one can check that $g(H) \subset H_{i,j}$ if and only if $H \cap
H_i \cap H_j$ is a codimension two flat in $\sL(\A)$ if and only if
$H_i \cap H_j \subset H$.  Using this observation, the relation
\eqref{eq:relation} may be expressed as
\[
[C^1_H,C^1_{H_m}] = 
[C^1_{H_m},\sum_{H_m \cap H_j \subset H} C^1_{H_j}].
\]
A calculation then shows that this relation is equivalent to
$[C^1_X,C^1_{H_m}]=0$, where $X$ is the codimension two flat in
$\sL(\A)$ contained in $H$ and $H_m$.  Since this relation holds for
all $H_m\in\A\setminus\B$ for which $X \subset H\cap H_m$, it follows
that $[C^1_X,C^1_{H}]=0$ as well.
\end{proof}

\section{Homology of the Loop Space} \label{sec:loop}

The structure of the Lie algebra of primitive elements in the homology
of the loop space of the complement of a redundant subspace
arrangement associated to a fiber-type hyperplane arrangement is
analyzed in this section.  In analogy with the previous section, begin
by recalling this structure for the classical configuration spaces
$F(\C^{k+1},n)$ for $k\ge 1$.

\begin{exm} \label{exm:config}
The integral homology of the loop space $\Omega F(\C^{k+1},n)$ was
calculated by Fadell and Husseini \cite{FH}.  The structure of the Lie
algebra $\Prim H_*(\Omega F(\C^{k+1},n);\Z)$ was subsequently
determined by Cohen and Gitler \cite{CG}.  For brevity, denote this
Lie algebra by $\cL(n)_{k}$.  The structure of $\cL(n)_{k}$ may be
described as follows.

For each $n \ge 2$, let $L[n]_k$ denote the free Lie algebra generated
by elements $B_{i,n+1}$, $1\le i \le n$, of degree $2k$.  Then
$\cL(n)_k$ is additively isomorphic to the direct sum
$\bigoplus_{j=1}^{n-1} L[j]_k$, and the Lie bracket relations in
$\cL(n)_k$ are given by the infinitesimal pure braid relations on the
$B_{i,j}$, see \eqref{eq:braidrels}.  Thus, there is an isomorphism of
graded Lie algebras $\cL(n)_{k} \cong E_0^*(P_n)_k$, see Definition
\ref{defn:EGk}, Theorem \ref{thm:config}, and Example \ref{exm:Pn}.

Furthermore, as is the case for the descending central series of the
pure braid group, the Lie algebra $\cL(n+1)_k$ is isomorphic to the
semidirect product of $\cL(n)_k$ by $L[n]_k$ determined by the Lie
homomorphism $\theta_n^k:\cL(n)_k \to \Der(L[n]_k)$ given by
$\theta_n^k(B_{i,j}) = \ad(B_{i,j})$.  From the infinitesmial pure
braid relations, there is a formula for $\ad(B_{i,j})$ analogous to
that given in \eqref{eq:adAij}.  As before, this extension of Lie
algebras arises topologically from the bundle of configuration spaces
$F(\C^{k+1},n+1)\to F(\C^{k+1},n)$.
\end{exm}

Now let $\A$ be a fiber-type hyperplane arrangement in $\C^\ell$ with
exponents $\{d_1,\dots,d_\ell\}$.  Then, for each $k$, there is a
tower of fiber bundles
\begin{equation*} \label{eq:BundleTower2}
M(\A^k_\ell) \xrightarrow{\ p^k_{\ell}\ } M(\A^k_{\ell-1})
\xrightarrow{p^k_{\ell-1}} \cdots \xrightarrow{\ p^k_{2}\ }
M(\A^k_1) = \C^k \setminus \{d_1 \text{ points}\},
\end{equation*}
where the fiber of $p^k_{j}$ is homeomorphic to the complement of
$d_{j}$ points in $\C^k$, see Theorem~\ref{thm:ftbundle}. 
Furthermore, each of the fiber bundles $p^k_j:M(\A^k_j) \to
M(\A^k_{j-1})$ involving the complements of the redundant subspace
arrangements $\A_j^k \subset \C^{jk}$ admits a cross-section, and, as
indicated above, $M(\A^k_1)$ is the complement of $d_1$ points in
$\C^k$.  By work of the second two authors \cite[Theorem 1]{CX}, the
following holds.

\begin{thm} \label{thm:loophomology}
Let $\A$ be a fiber-type hyperplane arrangement in $\C^\ell$ with
exponents $\{d_1,\dots,d_\ell\}$.  Then, for each $k \ge 1$,
\begin{enumerate}
\item[(a)] There is a homotopy equivalence $\Omega M(\A^{k+1}) \to
\prod_{j=1}^\ell \Omega(\C^{k+1} \setminus \{\text{$d_j$ points}\})$;

\smallskip

\item[(b)] The integral homology of $\Omega M(\A^{k+1})$ is torsion
free, and is isomorphic to the tensor product $\bigotimes_{j=1}^\ell
H_*(\Omega(\C^{k+1} \setminus \{\text{$d_j$ points}\});\Z)$ as a
coalgebra;

\smallskip

\item[(c)] The module of primitives in the integral homology of
$\Omega M(\A^{k+1})$ is isomorphic to
$\pi_*(\Omega M(\A^{k+1}))$ modulo torsion as a Lie algebra.
\end{enumerate}
\end{thm}

\begin{rem}
The homotopy groups of a loop space admit a bilinear pairing, which
satisfies the axioms for a graded Lie algebra in case there is no $2$
or $3$ torsion in the homotopy groups.  The graded analogue of the
symmetry law can fail in case $2$-torsion is present, while the graded
analogue of the Jacobi identity can fail if $3$-torsion is present. 
Thus, forming the quotient of the homotopy groups by the torsion gives
a graded module which satisfies the axioms for a graded Lie algebra. 
Analogous properties of iterated loop spaces yield a graded Poisson
algebra, see Section \ref{sec:poisson} below.
\end{rem}

\begin{proof}[Proof of Theorem \ref{thm:loophomology}]
Given a fibration $F \xrightarrow{\,i\,}E
\xrightarrow{\phantom{\,i\,}} B$ with a section $\sigma$, there is a
homotopy equivalence $\Omega E \simeq \Omega B \times \Omega F $ given
by the composite:
\[
\Omega B \times \Omega F \xrightarrow{\;\Omega \sigma \times \Omega i \;}
\Omega E \times \Omega E  \xrightarrow{\; \mu \; }  \Omega E,
\]
where $\mu$ is the loop space multiplication and such that the
inclusions of $\Omega B$ and $\Omega F$ into $\Omega E$ are maps of
$H$-spaces.  Morever, if the spaces involved have torsion free
homology then
$H_*(\Omega E) \cong H_*(\Omega B) \otimes H_*(\Omega F)$.  
By a theorem of Milnor and Moore, one obtains
\begin{equation}\label{eq:MilnorMoore}
\Prim H_*(\Omega E) \cong
            \Prim H_*(\Omega B) \oplus \Prim H_*(\Omega F)
\end{equation}
upon passing to the Lie algebra of primitives.  This result is a
topological analogue of Theorem~\ref{thm:FalkRandell} as the
underlying Lie algebra structure is ``twisted.''

Now apply these considerations to the fiber bundle $p^{k+1}_j :
M(\A^{k+1}_j) \to M(\A^{k+1}_{j-1})$.  The fiber in this case is $F =
\C^{k+1} \setminus \{\text{$d_j$ points} \} \simeq \bigvee_{d_j}
S^{2k+1}$.  Assertion (a) follows by induction, and then (b) by the
K\"unneth theorem.  By the Bott-Samelson Theorem, $H_*(\Omega F)$ is
isomorphic to $T[d_j]_k$, a tensor algebra on $d_j$ generators of
degree $2k$.

Thus the module of primitive elements is generated as a Lie algebra by
the primitive elements in degree $2k$ which are in the image of the
Hurewicz map.  Since the Hurewicz map takes values in the the module
of primitive elements, that module is generated as a Lie algebra by
those spherical classes given by the homology classes of degree $2k$.

Next notice that the homology groups here are torsion free.  Hence the
Hurewicz map factors through $\pi_*\Omega(M(\A^{k+1}))/\Tors$. 
Furthermore, the homotopy groups of a loop space modulo torsion give a
graded Lie algebra where the Lie bracket is induced by the classical
Samelson product, and the Hurewicz map is a morphism of graded Lie
algebras.  Thus the induced map $\pi_*\Omega(M(\A^{k+1}))/\Tors \to
\Prim H_*(\Omega M(\A^{k+1});\Z)$ is an epimorphism of Lie algebras.

Since all spaces are simply connected, and are of finite type, the
homotopy groups modulo torsion are finitely generated free abelian
groups in any fixed degree.  By a classical theorem of Milnor and
Moore concerning rational homotopy groups, the induced map
$\pi_*\Omega(M(\A^{k+1}))/\Tors \to \Prim H_*(\Omega M(\A^{k+1});\Z)$
is also a monomorphism.  The result follows.
\end{proof}

By \eqref{eq:MilnorMoore} above, there is an isomorphism of graded
abelian groups
\[
\Prim H_*(\Omega M(\A^{k+1}_j))  \cong
\Prim H_*(\Omega M(\A^{k+1}_{j-1})) \oplus  L[d_j]_k.
\]
Proceeding inductively, this implies that the Lie algebra $\Prim
H_*(\Omega M(\A^{k+1});\Z)$ is isomorphic to $L[d_1]_k\oplus \cdots
\oplus L[d_\ell]_k$ as a graded abelian group, where
$\{d_1,\dots,d_\ell\}$ are the exponents of $\A$.  Thus the Lie
algebras $\Prim H_*(\Omega M(\A^{k+1});\Z)$ and $E_0^*(G(\A))_k$ are
additively isomorphic, see Theorem \ref{thm:DCSadditive}.  To show
that they are are isomorphic as Lie algebras, thereby completing the
proof of Theorem \ref{thm:CXconj}, it remains to show that the Lie
bracket structure of $\Prim H_*(\Omega M(\A^{k+1});\Z)$ coincides with
that of $E_0^*(G(\A))_k$.  This analysis parallels the determination
of the Lie algebra structure of $E_0^*(G(\A))$ in
Section~\ref{sec:DCS}.

The fiber-type hyperplane arrangement $\A=\A_{\ell}$ is strictly
linearly fibered over $\A_{\ell-1}$, and
$|\A|=|\A_{\ell-1}|+d_{\ell}$.  As before, write $\B=\A_{\ell-1}$ and
$n=d_{\ell}$.  Recall the map $g^{k+1}:M(\B^{k+1})\to F(\C^{k+1},n)$
from \eqref{eq:kroot}.  Recall also that the Lie algebra $\Prim
H_*(\Omega F(\C^{k+1},n);\Z)$ is denoted by $\cL(n)_k$.  Analogously,
denote the Lie algebra $\Prim H_*(\Omega M(\A^{k+1});\Z)$ by
$\cL(\A)_k$.

\begin{thm} \label{thm:loopstructure}
Let $\A$ and $\B$ be fiber-type hyperplane arrangements with $\A$
strictly linearly fibered over $\B$ and $|\A|=|\B|+n$.  Then the Lie
algebra $\cL(\A)_k$ is isomorphic to the semidirect product of
$\cL(\B)_k$ by the free Lie algebra $L[n]_k$ determined by the Lie
homomorphism $\Theta^{k} = \theta_{n}^{k} \circ
\c_{*}^{k}:\cL(\B)_k\to \Der(L[n]_{k})$, where
$\c_{*}^{k}:\cL(\B)_k\to \cL(n)_k$ is the map in loop space homology
induced by $g^{k+1}:M(\B^{k+1})\to F(\C^{k+1},n)$, and
$\theta_n^{k}:\cL(n)_k \to \Der(L[n]_{k})$ is given by
$\theta_n^{k}(B_{i,j}) = \ad(B_{i,j})$.
\end{thm}
\begin{proof}
The realization of the bundle $p^{k+1}:M(\A^{k+1})\to M(\B^{k+1})$ as
the pullback of the bundle of configuration spaces
$p_{n+1}^{k+1}:F(\C^{k+1},n+1)\to F(\C^{k+1},n)$ along the map
$g^{k+1}:M(\B^{k+1})\to F(\C^{k+1},n)$ from Theorem
\ref{thm:slfsoupup} yields a commutative diagram of Hopf algebras
\[
\begin{CD}
H_{*}(\Omega(\C^{k+1}\setminus \{n\ \text{points}\}))    @>>>
H_{*}(\Omega M(\A^{k+1})) @>>> H_{*}(\Omega M(\B^{k+1})) \\
@VV{\id}V      @VVV  @VV{\c_{*}^{k}}V\\
H_{*}(\Omega(\C^{k+1}\setminus \{n\ \text{points}\}))
@>>> H_{*}(\Omega F(\C^{k+1},n+1)) @>>>H_{*}(\Omega F(\C^{k+1},n))
\end{CD}
\]
with exact rows, and, on the level of primitives, a commutative
diagram of Lie algebras
\[
\begin{CD}
0 @>>> L[n]_{k}    @>>> \cL(\A)_k @>>> \cL(\B)_k @>>>0\\
@.     @VV{\id}V      @VVV  @VV{\c_{*}^{k}}V\\
0 @>>> L[n]_{k}    @>>> \cL(n+1)_k @>>> \cL(n)_k @>>>0
\end{CD}
\]
where $\c_{*}^{k}:\cL(\B)_k \to \cL(n)_k$ is induced by
$g^{k+1}:M(\B^{k+1})\to F(\C^{k+1},n)$.  Since the underlying bundles
admit cross-sections, the rows in the above diagrams are split exact.
Via these splittings, view $\cL(\B)_k$ and $\cL(n)_k$ as Lie
subalgebras of $\cL(\A)_k$ and $\cL(n+1)_k$ respectively.

 From the above considerations, it follows that the Lie algebra
$\cL(\A)_{k}$ is isomorphic to the semidirect product of $\cL(\B)_{k}$
by $L[n]_{k}$ determined by the Lie homomorphism
$\Theta^{k}:\cL(\B)_{k} \to \Der(L[n]_{k})$ given by
$\Theta^{k}(b)=\ad_{L[n]_{k}}(b)$ for $b \in \cL(\B)_{k}$.  Moreover,
for $a\in L[n]_{k}$, we have $[b,a]=[\c^{k}_*(b),a]$ in $L[n]_{k}$.
Thus $\ad_{L[n]_{k}}(b) = \ad_{L[n]_{k}}(\c^{k}_*(b))$ in
$\Der(L[n]_{k})$ and $\Theta^{k}=\theta^{k}_n\circ \c^{k}_*$.
\end{proof}

This result, together with Proposition \ref{prop:inducedmap}, provides
an inductive description of the Lie bracket structure of $\cL(\A)_{k}$. 
The space $M(\B^{k+1})$ is $2k$-connected, and the cohomology algebra
$H^*(M(\B^{k+1});\Z)$ is generated by classes $a^{k+1}_H$ in
one-to-one correspondence with the hyperplanes $H\in \B$, see
Corollary \ref{cor:cohomology}.  These classes are of degree $2k+1$,
and are dual to the elements of the basis $\{C^{k+1}_H \mid H\in \B\}$
for $H_{2k+1}(M(\B^{k+1});\Z)$ exhibited in Proposition
\ref{prop:basis}.  See also Remark \ref{rem:dual}.

The above observations imply that homology suspension induces an
isomorphism
\[
\sigma_*:H_{2k}(\Omega M(\B^{k+1});\Z) \to H_{2k+1}(M(\B^{k+1});\Z).
\]
Let $\beta^{k}_H \in H_{2k}(\Omega M(\B^{k+1});\Z)$ be the unique
class satisfying $\sigma_*(\beta^{k}_H)=C^{k+1}_H$.  Recall that the
free Lie algebra $L[n]_{k}$ is generated by
$B_{1,n+1},\dots,B_{n,n+1}$.

\begin{cor} \label{cor:kbrackets}
For generators $\beta^{k}_H$ of $\cL(\B)_{k}$ and $B_{m,n+1}$ of
$L[n]_{k}$, one has
\[
\Theta^{k}(\beta^{k}_H)(B_{m,n+1}) = \sum_{g(H) \subset H_{i,j}}
[B_{i,j},B_{m,n+1}].
\]
\end{cor}
\begin{proof}
By Proposition \ref{prop:inducedmap}, one has
$g^{k+1}_{*}(C^{k+1}_{H})= \sum A_{i,j}$, where the sum is over all
$i$ and $j$ for which $g(H) \subset H_{i,j}$.  Since the homology
suspension $\sigma_{*}$ is an isomorphism and $\c_{*}^{k}$ is the map
in loop space homology induced by $g^{k+1}$, one has 
$\c_{*}^{k}(\beta^{k}_{H})= \sum B_{i,j}$, where the sum is over all
$i$ and $j$ for which $g(H) \subset H_{i,j}$.  The result follows.
\end{proof}

To complete the proof of Theorem \ref{thm:CXconj}, assume inductively
that the Lie algebras $E^{*}_{0}(G(\B))_{k}$ and $\cL(\B)_{k}$ are
isomorphic.  By Theorem \ref{thm:DCSstructure}, the Lie algebra
$E^{*}_{0}(G(\A))$ is the extension of $E^{*}_{0}(G(\B))$ by the free
Lie algebra $L[n]$ (generated in degree one) determined by the Lie
homomorphism $\Theta=\theta_{n} \circ g_{*}$.  Thus
$E^{*}_{0}(G(\A))_{k}$ may be realized as the extension of
$E^{*}_{0}(G(\B))_{k}$ by the free Lie algebra $L[n]_{k}$ (generated
in degree $2k$) determined by $\Theta$ as specified in Definition
\ref{defn:EGk}.  Similarly, by Theorem \ref{thm:loopstructure}, the
Lie algebra $\cL(\A)_{k}$ is the extension of $\cL(\B)_{k}$ by the
free Lie algebra $L[n]_{k}$ determined by the Lie homomorphism
$\Theta^{k}=\theta_{n}^{k} \circ \c^{k}_{*}$.  A comparison of the
results of Corollary~\ref{cor:brackets} and
Corollary~\ref{cor:kbrackets} reveals that these extensions coincide. 
Therefore, the Lie algebras $E^{*}_{0}(G(\A))_{k}$ and $\cL(\A)_{k}$
are isomorphic.

Alternatively, Corollary \ref{cor:kbrackets} may be used to explicitly
determine the Lie bracket structure in $\cL(\A)_k$.  As the argument
is completely analogous to that which established
Theorem~\ref{thm:LieBrackets}, the result stated below without proof. 
The Lie algebra $\cL(\A)_k$ is generated by $\{\beta^k_H \mid
H\in\A\}$.  For a flat $X\in \sL(\A)$, write
$\beta^{k}_{X}=\sum_{X\subset H} \beta^k_H$.

\begin{thm} \label{thm:kLieBrackets}
Let $\A$ be a fiber-type hyperplane arrangement.  Then, for each $k\ge
1$, the Lie bracket relations in $\cL(\A)_k$ are given by
\[
[\beta^k_X,\beta^k_H]=0,
\]
for codimension two flats $X\in\sL(\A)$ and hyerplanes $H\in\A$ 
containing $X$.
\end{thm}

\section{Homology of Iterated Loop Spaces} \label{sec:poisson}

In this final section, the Poisson algebra structure on the homology
of an iterated loop space of the complement of a redundant subspace
arrangement associated to a fiber-type hyperplane arrangement is
briefly analyzed.

For $q>1$, the homology of an $q$-fold loop space, $\Omega^q X$,
admits the structure of a graded Poisson algebra.  Namely, there is a
bilinear map given by the Browder operation
\[
\lambda_{q-1}: H_i(\Omega^q X) \otimes H_j(\Omega^q X) \to
H_{i+j+q-1} (\Omega^q X)
\]
which satisfies properties listed in \cite[pages 215--217]{CLM}.  In
particular, this pairing satisfies the axioms of a (graded) Poisson
algebra, and is compatible with the Whitehead product structure for
the classical Hurewicz homomorphism.

In the case where $X=X_\ell(\C^{k+1})$, $k\ge 1$, satisfies conditions
(1)--(3) from the Introduction, these structures are analogues of
classical constructions in homotopy theory.  First, note that the
single suspension $\Sigma X_\ell(\C^{k+1})$ is homotopy equivalent to
a bouquet of spheres.  Thus there is an induced map $\sigma^2:
\Sigma^2 X_\ell(\C^{k+1}) \to X_\ell(\C^{k+2})$ which induces an
isomorphism on the first non-trivial homology group.  The adjoint
$E^2: X_\ell(\C^{k+1}) \to\ \Omega^2 X_\ell(\C^{k+2})$ also induces an
isomorphism on the first non-trivial homology group.  This last map is
an analogue of the classical Freudenthal double suspension where the
spaces $X_\ell(\bE)$ are replaced by single odd dimensional spheres. 
Looping $E^2$ is given by $\Omega(E^2): \Omega X_\ell(\C^{k+1}) \to\
\Omega^3 X_\ell(\C^{k+2})$.

\begin{thm} \label{thm:ftpoisson}
Let $\A$ be a fiber-type hyperplane arrangement in $\C^\ell$ with
exponents $\{d_1,\dots,d_\ell\}$.  Then, for each $k \ge 1$,
\begin{enumerate}
\item[(a)] The multiplicative map $\Omega(E^2): \Omega M(\A^{k+1}) \to
\Omega^3 M(\A^{k+2})$ induces an isomorphism on $H_{2k}(-;\Z)$, and is
zero in degrees greater than $2k$.

\smallskip

\item[(b)] If $q > 1$, the homology of $\Omega^q M(\A^{k+1})$, with any
field coefficients, is a graded Poisson algebra with Poisson bracket
given by the Browder operation for the homology of a $q$-fold loop
space.

\smallskip

\item[(c)] If $q>1$, then $\Omega^q M(\A^{k+1})$ is homotopy
equivalent to $\prod_{j=1}^{\ell} \Omega^q (\bigvee_{d_{j}}
S^{2k+1})$.

\smallskip

\item[(d)] If $1<q<2k+1$, the homology of $\Omega^{q} M(\A^{k+1})$,
with coefficients in a field $\bF$ of characteristic zero, is
generated as a Poisson algebra by elements $\beta_{H}$ of degree
$2k+1-q$ for $H\in\A$.  The Poisson bracket is given by the Browder
operation $\lambda_{q-1}$, and satisfies the relations
\[
\lambda_{q-1}[\beta_X,\beta_H],
\]
for codimension two flats $X\in\sL(\A)$ and hyerplanes $H\in\A$
containing $X$, where $\beta_X = \sum_{X\subset H} \beta_H$.
\end{enumerate}
\end{thm}

\begin{proof}[Sketch of Proof]
Part (a) follows from the fact that the homology of $\Omega^3
M(\A^{k+2})$ is abelian while the homology of $\Omega M(\A^{k+1})$ is
generated by Lie brackets of weight at least $2$ in homological
degrees greater than $2k$.

Part (b) follows from the remarks at the beginning of this section.

Part (c) follows at once from the fact that the result holds in case
$q = 1$, which was established in Theorem \ref{thm:loophomology}.

In case $ q = 1 $, the Browder operation $\lambda_{q-1}$ is precisely
the natural Lie bracket in the homology of a $1$-fold loop space,
$\Omega M(\A^{k+1})$.  These Lie bracket relations are recorded in
Theorem~\ref{thm:kLieBrackets}.  As shown in \cite[pages
215--217]{CLM}, a further property of the operation $\lambda_{q-1}$ is
that $ \sigma_* \lambda_{q-1}(x,y) = \lambda_{q-2}(\sigma_* x,\sigma_*
y)$, where $\sigma_*$ denotes the homology suspension.  Thus by
induction on $q$, the asserted Poisson bracket relations are satisfied
modulo elements in the kernel of the suspension.  Furthermore,
$\lambda_{q-1}(x,y)$ is primitive in case the classes $x$ and $y$ are
primitive.

In characteristic zero, and in case $q$ is greater than $1$, the
homology suspension induces an isomorphism on the module of
primitives.  Thus the asserted Poisson bracket relations are
satisfied.
\end{proof}

\begin{rem} Let $\A=\A_{n}$ be the braid arrangement in $\C^{n}$.  As
noted in Example~\ref{exm:braidarrangement}, one then has
$M(\A_{n}^{k+1})=F(\C^{k+1},n)$ for all $k$.  For the braid
arrangement, the codimension two flats in $\sL(\A_{n})$ (the partition
lattice) are of the forms
\[
H_{i,j}\cap H_{i,k}\cap H_{j,k}\ \text{for}\ 1\le i<j<k\le n, \quad
\text{and}\quad H_{i,j}\cap H_{k,l}\ \text{for}\
\{i,j\}\cap\{k,l\}=\emptyset.
\]
Thus by Theorem \ref{thm:ftpoisson}, for $1<q<2k+1$, the homology of
$\Omega^{q} F(\C^{k+1},n)$, with coefficients in a field $\bF$ of
characteristic zero, is generated as a Poisson algebra by elements
$B_{i,j}=\beta_{H_{i,j}}$ of degree $2k+1-q$ for $1\le i<j \le n$. 
Moreover, the Poisson bracket relations are given by the universal
infinitesimal Poisson braid relations:
\begin{equation*} 
\begin{aligned}
\lambda_{q-1}{[}B_{i,j}+B_{i,k}+B_{j,k},B_{m,k}]&=0 \quad
\text{for $m=i$ or $m=j$, and}\\
\lambda_{q-1}[B_{i,j},B_{k,l}]&=0
\quad\text{for $\{i,j\}\cap\{k,l\}=\emptyset$.}
\end{aligned}
\end{equation*}
As shown in \cite{CG}, these are precisely the infinitesimal pure
braid relations in case $q=1$, see also Examples \ref{exm:Pn} and
\ref{exm:config}.

It seems likely that, via the natural universal mapping property, one
could define the ``universal infinitesimal Poisson braid algebra,''
and that the homology of $\Omega^q F(\C^{k+1},n)$ with coefficients in
a field $\bF$ of characteristic zero is that algebra over $\bF$.
\end{rem}

\bibliographystyle{amsalpha}

\end{document}